# A severe inconsistency of transfinite set theory[1]

W. Mückenheim

University of Applied Sciences, D-86161 Augsburg, Germany [mueckenh@rz.fh-augsburg.de]

Transfinite set theory including the axiom of choice supplies the following basic theorems: (1) Mappings between infinite sets can always be completed, such that at least one of the sets is exhausted. (2) The real numbers can be well ordered. (3) The relative positions of real numbers which are enumerated by natural numbers can always be determined, in particular the maximum real number below a given limit. (4) Any two different real numbers are separated by at least one rational number.
These theorems are applied to map the irrational numbers into the rational numbers, showing that the set of all irrational numbers is countable.

## 1. Introduction

Georg Cantor is that mathematician whose name is most closely connected with the notion of infinity. An essential point of his theory is the basic law (Denkgesetz) that the continuum can be well-ordered. Precisely 100 years ago, here in Heidelberg, Cantor had to suffer a serious attack from König[2] who claimed to have contradicted this. In fact nobody has ever accomplished it. But Zermelo's axiom of choice, created immediately after that dispute, convinced more and more mathematicians that well-ordering could be proved. This is a singular case in history of mathematics. Instead of refuting König's claim by a counter example, an axiom was created which supplies rather an encoded proclamation of the desired result than proof[3].

The axiom of choice allows us to manipulate all the elements of an infinite set simultaneously instead of considering one after the other[4]. Any mapping between two infinite sets leads to the exhaustion of at least one of these sets. So all infinite processes can be executed completely, or, to put it in other, perhaps heretical words, all infinities can be finished.

In the following we will see that this property of the axiom of choice leads to a severe inconsistency of transfinite set theory, contradicting its most important theorem according to which there were more irrational numbers than rational numbers.

---

[1] Talk delivered at the meeting of the German Math. Soc., section logic, in Heidelberg, Sept. 14, 2004.

[2] In his talk at the third international congress in 1904 König had used Bernstein's theorem $\aleph_x^{\aleph_0} = \aleph_x 2^{\aleph_0}$ for $x = \omega$. But it holds only in the domain of finite $x$, as König admitted in 1905 [König a]. In that very year, however, he published further evidence against well-ordering: Die bisher entwickelten Annahmen führen in merkwürdig einfacher Weise zu dem Schlusse, daß *das Kontinuum nicht wohlgeordnet werden kann. ...* Die Annahme, daß das Kontinuum wohlgeordnet werden kann, hat demnach zu einem Widerspruch geführt. [König b].

[3] Zermelo claimed in the head-line of his paper to have proven that any set can be well-ordered: Beweis, daß jede Menge wohlgeordnet werden kann [Zermelo].

[4] Zur Unterstützung dieser zeitlosen Auffassung hat E. Zermelo den glücklichen Gedanken gehabt, von vornherein aus jeder von Null verschiedenen Teilmenge A' von A eins ihrer Elemente a' = f(A') auszuwählen ... Das System sukzessiver Wahlakte ist damit durch ein, in der Praxis des Denkens natürlich ebenso unausführbares System simultaner Wahlakte ersetzt. [Hausdorff, p. 133]

## 2. Preparation

The rational numbers are countable whereas the irrational numbers are uncountable. It is argued that they make up intervals of the continuum separated by rational numbers which are only points. The set of intervals formed by these points necessarily is countable too. In order to support the idea of uncountably many irrational numbers, we should be able to find at least one interval containing uncountably many of them. This task, however, is condemned to fail, because between any two different real numbers $r_1 < r_2$ at least one rational number $q$ can be found separating them. Otherwise $r_1$ and $r_2$ would not be different but one and the same number. This is easy to prove. Moreover it is a fundamental part of the set-theoretic definition of the order type of the real numbers. By human sense or internal vision one could conclude then that at least half of the real numbers must be rational. Though Cantor himself frequently referred to internal vision[5] this is not a decisive argument. Therefore I will give firm proof. In order to use only correct arguments I will collect them from the following proof given by Cantor himself.

Let $\{\Delta_\nu\} \subset (0, 1)$ be an infinite but countable set of disjoint intervals, ordered by magnitude, and let $\{\varphi_\nu\}$ be any dense (überalldicht) sequence of different points $\{\varphi_\nu\} \subset (0, 1)$. Between $\{\varphi_\nu\}$ and $\{\Delta_\nu\}$ a one-to-one correspondence or bijection can be established, such that the relative position (Lagenbeziehung) of $\varphi_{\kappa_\nu}$ and $\varphi_{\kappa_\mu}$ is the same as that of $\Delta_\nu$ and $\Delta_\mu$.[6] This bijection can be accomplished as follows by the mapping: $\varphi_{\kappa_1} \to \Delta_1$ where $\varphi_{\kappa_1}$ is the first point $\varphi_1$ of $\{\varphi_\nu\}$. $\varphi_{\kappa_2} \to \Delta_2$ where $\varphi_{\kappa_2}$ is the first point of $\{\varphi_\nu\}$ which has the same relative position (Lagenbeziehung) with respect to $\varphi_{\kappa_1}$ as has $\Delta_2$ with respect to $\Delta_1$. ... $\varphi_{\kappa_{\nu+1}} \to \Delta_{\nu+1}$ where $\varphi_{\kappa_{\nu+1}}$ is the first point of $\{\varphi_\nu\}$ which has the same relative position (Lagenbeziehung) with respect to $\varphi_{\kappa_1}, \varphi_{\kappa_2}, ..., \varphi_{\kappa_\nu}$ as has $\Delta_{\nu+1}$ with respect to $\Delta_1, \Delta_2, ..., \Delta_\nu$.

This mapping can be completed such that the whole infinite set is included. In this way *all intervals* $\Delta_\nu$ are related to definite points $\varphi_\nu$.[7] The proof is continued by showing that this is really a one-to-one correspondence, but we do not need that property.

---

[5] ... und da andererseits der Anzahlbegriff in unserer inneren Anschauung eine unmittelbare gegenständliche Repräsentation erhält ... [Cantor, p. 168] ... wie die Gesetze derselben aus der unmittelbaren inneren Anschauung mit apodiktischer Gewißheit erschlossen werden. [Cantor, p. 170] ... wie man aus der inneren Anschauung unmittelbar erkennt. [Cantor, p. 201]

[6] Zwischen der Punktmenge $\{\varphi_\nu\}$ einerseits und der Intervallmenge $\{\Delta_\nu\}$ andrerseits läßt sich immer eine solche gesetzmäßige, gegenseitig eindeutige und vollständige Korrespondenz ihrer Elemente herstellen, daß, wenn $(\Delta_\nu)$ und $(\Delta_\mu)$ irgend zwei Intervalle der Intervallmenge, $\varphi_{\kappa_\nu}$ und $\varphi_{\kappa_\mu}$ die zu ihnen gehörigen Punkte der Punktmenge $\{\varphi_\nu\}$ sind, alsdann stets $\varphi_{\kappa_\nu}$ links oder rechts von $\varphi_{\kappa_\mu}$ liegt, je nachdem das Intervall $(\Delta_\nu)$ links oder rechts von dem Intervalle $(\Delta_\mu)$ fällt, oder, was dasselbe heißen soll, daß die Lagenbeziehung der Punkte $\varphi_{\kappa_\nu}$ und $\varphi_{\kappa_\mu}$ stets dieselbe ist wie die Lagenbeziehung der ihnen entsprechenden Intervalle $(\Delta_\nu)$ und $(\Delta_\mu)$. [Cantor, p. 238] In Cantor's original text the interval $\Delta_\nu$ is denoted by $(a_\nu...b_\nu)$.

[7] Klar ist zunächst, daß auf diese Weise *allen* Intervallen ... *bestimmte* Punkte ... zugeordnet werden; denn wegen des Überalldichtseins der Menge $\{\varphi_\nu\}$ im Intervall $(0 ... 1)$, und weil die Endpunkte 0 und 1 nicht zu $\{\varphi_\nu\}$ gehören, gibt es in dieser Reihe unendlich viele Punkte, die eine *geforderte* Lagenbeziehung zu einer bestimmten endlichen Anzahl von Punkten derselben Menge $\{\varphi_\nu\}$ besitzen, und es erfährt daher der aus unsrer Regel resultierende Zuordnungsprozeß *keinen Stillstand*. [Cantor, p. 239]

From Cantor's proof we obtain as a basic theorem that the relative positions (Lagenbeziehung) of real numbers in well-ordered countable sets can always be determined. This is due to the fact that any element of a countable set can be enumerated by a natural number, and each natural number is finite. So we always have to deal with finite indices. Otherwise set theory would break down. As an example consider a Cantor-list of real numbers enumerated by the set $\mathbb{N}$ of natural numbers $n$. For an infinite line number $\omega$, we would not be able to determine the corresponding diagonal element $a_{\omega\omega}$ of the Cantor-list to be exchanged. Hence any number $n$ enumerating a line of the Cantor-list is finite and can be distinguished from its neighbours.[8]

## 3. Theorem and proof

The basic theorems obtained from the foregoing can be summarised as follows:

- Mappings between infinite sets can always be completed, such that at least one of the sets is exhausted.
- The real numbers can be well ordered.
- The relative positions (Lagenbeziehung) of real numbers enumerated by natural numbers can always be determined, in particular the maximum real number below a given limit.
- Between any two real numbers, there exists always a rational number.

They are sufficient to prove that the set of all irrational numbers can be mapped *into* the set of all rational numbers, showing that the cardinality of the former does *not* surpass that of the latter. For the sake of simplicity I will restrict this proof to positive numbers. The extension to all numbers is then obvious.

**Theorem**. *The set of all positive irrational numbers $\mathbb{X}_+$ can be mapped into the set of all positive rational numbers $\mathbb{Q}_+$ leading to $|\mathbb{X}_+| \leq |\mathbb{Q}_+|$.*

**Proof.** Let the sets $\mathbb{Q}_+$ and $\mathbb{X}_+$ be well-ordered. Define two sets, one of them containing only the number zero, $Q = \{0\}$, and the other one being empty, $X = \{\}$. Take the first element $\xi_1 \in \mathbb{X}_+$. Select the largest rational number $q \in Q$ with $q < \xi_1$ (in the first step, this is obviously $q = 0$). Between two different real numbers like $q$ and $\xi_1$ there is always a rational number. There are many of them. Take the first one $q_1 \in \mathbb{Q}_+$, with $q < q_1 < \xi_1$. Transfer $\xi_1$ to the set $X$ and transfer $q_1$ to the set $Q$. Then choose the next positive irrational number, $\xi_2 \in \mathbb{X}_+$, select the largest $q \in Q = \{0, q_1\}$ with $q < \xi_2$. Take the first rational number, $q_2 \in \mathbb{Q}_+$, with $q < q_2 < \xi_2$. Transfer $\xi_2$ to the set $X$ and transfer $q_2$ to the set $Q$. Continue until one of the sets $\mathbb{Q}_+$ or $\mathbb{X}_+$ is exhausted which, according to the axiom of choice, will unavoidably occur.

If the set $\mathbb{Q}_+$ were exhausted prematurely and no $q_n$ remained available to map $\xi_n$ on it, this proof would fail. We would leave the countable domain and could no longer make use of Cantor's

---

[8] By the way, here lies the reason why Cantor's diagonal argument must fail, nevertheless. Every number smaller than $\omega$ is a finite number, and it is surpassed by other finite numbers [Cantor, p. 406]. Hence, if the list contains "*every* number smaller than $\omega$", then there are *other* numbers, not contained in the list. The list is not complete. Further, for finite numbers the miracle of the discrepancy between ordinal and cardinal number is void of its witchcraft. Hence the cardinal number of the lines of the list is finite, however large the finite ordinals may be. Therefore the diagonal number constructed is not a transfinite number, not even an irrational number, but it is simply a rational number, because of its finite number of digits.

"Lagenbeziehung" to select the largest rational number $q \in Q$ with $q < \xi_n$. But that cannot occur because there is *always* a rational number between two real numbers. As long as rationals $q_n \in Q_+$ are available, the set of pairs $(q_n, \xi_n)$ remains countable and there are also natural numbers $n$ available as indices, because all our positive rational numbers have been enumerated by natural numbers. Therefore, we do not leave the countable domain and do not need transfinite induction. But our mapping process runs until one of the sets is exhausted. By tertium non datur this set can only be $X_+$.

The mapping supplies $|Q| = |X|$ at every stage while finally $X_{\text{fin}} = X_+$ and $Q_{\text{fin}} \subseteq Q_+$. ($Q_{\text{fin}}$ even must be a true subset of $Q_+$ because there are further rational numbers between every two elements of $Q_{\text{fin}}$.) The result $|X_+| \leq |Q_+|$ completes the proof.

## 4. Conclusion

Concluding we can say that there are no different infinities. If the axiom of choice is abolished, then well-ordering of the continuum and of larger sets is impossible, and there is no chance of attributing a cardinal number to those sets. If the axiom of choice is maintained then the continuum can be proved countable, also contradicting transfinite set theory.

Infinity marks a direction and not an actual magnitude or value. Thus, we need no alephs and no omegas and the continuum hypothesis is not only undecidable but completely meaningless. We need only one symbol to represent infinity, namely that one introduced by Wallis in 1655, and we can write
$$\infty = \infty + 1 = 2\infty = 2^{\infty}$$

and many more equations of this kind, which are tolerated by a direction.